\documentclass{amsart}
\usepackage{amscd,amssymb}

\theoremstyle{plain}
\newtheorem{thm}{Theorem}
\newtheorem{lem}[thm]{Lemma}

\newtheorem{prop}[thm]{Proposition}

\theoremstyle{definition}

\theoremstyle{remark}

\newtheorem{Remark}[thm]{Remark}

\numberwithin{equation}{section}

\begin{document} 

\title{ Complexity of countable categoricity in finite languages} 
\author{ Aleksander Ivanov  } 

\maketitle

\begin{quote} 
{\bf Abstract.} 
We study complexity of the index set of countably categorical 
theories and Ehrenfeucht theories in finite languages. 
\end{quote} 

S.Lempp and T.Slaman proved in \cite{LS} that indexes of 
decidable $\omega$-categorical theories form a $\Pi^0_3$-subset 
of the set of indexes of all computably enumerable theories. 
Moreover there is an infinite language so that the property 
of $\omega$-categoricity distinguishes a $\Pi^0_3$-complete 
subset of the set of indexes of computably enumerable 
theories of this language. 
Steffen Lempp asked the author if this could be done in a finite language. 
In this paper we give a positive answer (see Section 4). 
The crucial element of our proof is a theorem of Hrushovski on 
coding of $\omega$-categorical theories in finite languages 
(see \cite{hodges}, Section 7.4, pp. 353 - 355). 
Since we apply the method which was used in the the proof of 
this theorem, we present all the details in Section 1.  
Sections 2 - 3 contain several other applications of this theorem. 
In particular in the very short Section 2 we give an example 
of a non-G-compact $\omega$-categorical theory in a finite language. 
In Section 3 we show that there is a finite language such that 
the indexes of Ehrenfeucht theories with exactly three countable 
models form a $\Pi^1_1$-hard set. 
Here we also use the idea of Section 4 of \cite{LS} where a 
similar statement is proved in the case of infinite languages.  

The main results of the paper are available both 
for computability theorists and model theorists. 
The only place where a slightly advanced model-theoretical 
material appears is Section 2. 
On the other hand the argument applied in this section 
is very easy and all necessary preliminaries are presented.

\section{Hrushovski on $\omega$-categorical structures and finite languages} 

The material of this section is based on Section 7.4 of 
\cite{hodges}, pp. 353 - 355 (and preliminary notes of W.Hodges). 
We also give some additional modifications and remarks. 

Let $N$ be a structure in the language $L$ with a unary predicate $P$. 
For any family of relations $\mathcal{R}$ on $P$ definable in $N$ 
over $\emptyset$ one may consider the structure $M=(P,\mathcal{R})$. 
We say that $M$ is a {\em dense relativised reduct} if the image of 
the homomorphism $Aut(N)\rightarrow Aut(M)$ (defined by restriction) 
is dense in $Aut(M)$.  
\parskip0pt

Let $L$ be the language consisting of four unary symbols 
$P,Q,\lambda ,\rho$, a two-ary symbol $H$ and a four-ary one $S$. 
We will consider only $L$-structures where $P$ and $Q$ define a partition 
of the basic sort and $\lambda$, $\rho$ and $H$ are defined on $Q$. 
Moreover when $S(a,b,c,d)$ holds we have that $a,c\in P$ and $b,d \in Q$. 

\begin{thm} \label{Udi} 
If $M_0$ is any countable $\omega$-categorical structure then there 
is a countable $\omega$-categorical $L$-structure $N$ such that $M_0$ 
is a dense relativised reduct of $N$.  
In particular $M_0$ is interpretable in $N$ over $\emptyset$. 

For every set of sentences $\Phi$ axiomatising $Th(M_0 )$ the theory 
$Th(N)$ is axiomatised by a set of axioms which is computable with 
respect to $\Phi$ and the Ryll-Nardzewski function of $Th(M_0 )$.   
\end{thm} 

{\em Proof} (E.Hrushovski). 
Let $M_0$ be any countable $\omega$-categorical 
structure in a language $L_0$. 
We remind the reader that the {\em Ryll-Nardzewski function}  
of an $\omega$-categorical theory $T$ assigns 
to any natural $n$ the number of $n$-types of $T$. 
So by the set $\Phi$ as in the formulation and by the Ryll-Nardzewski 
function of $Th(M_0 )$ one can find an effective list of all pairwise 
non-equivalent formulas.  
Thus w.l.o.g. we may assume that $L_0$ is 1-sorted, relational and 
$M_0$ has quantifier elimination. 
In fact we can suppose that $L_0 = \{ R_1 ,R_2 ,..., R_n , ...\}$ 
where each $R_n$ describes a complete type in $M_0$ of arity 
not greater than $n$.  
We may also assume that for $m<n$ the arity of $R_m$ is not greater 
than the arity of $R_n$. 
We admit that tuples realising $R_n$ may have repeated coordinates. 

We now use standard material about Fra\"{i}ss\'{e} limits, 
see \cite{evans}. 
Note that the class of all finite substructures of $M_0$ (say $\mathcal{K}_0$) 
has the joint embedding and the amalgamation properties. 
Moreover for every $n$ the number of finite substructures of 
size $n$ is finite (this is the place where we use the assumption 
that each $R_n$ describes a complete type). 

Let us consider structures of the language $L\cup L_0$ which 
satisfy the property that all the relations $R_n$ are defined on $P$. 
For such a structure $M$ we call a tuple $(a_0 ,...,a_{m-1} ,c_0 ,...,c_{n-1} )$ 
of elements of $M$, an $n$-{\em pair of arity} $m$ if : 

(1) $m\le n$ and $M\models \bigwedge \{ P(a_i ): i<m\} \wedge \bigwedge \{ Q(c_j ): j<n\}$; 

(2) the elements $c_i$ are paiwise distinct and $M\models H(c_i, c_j )$ iff $(j=i+1) mod(n)$;  

(3) $M\models \lambda (c_i )$ iff $i =0$ and $M\models \rho (c_i )$ iff $i =m-1$; 

(4) $M\models S(a_i ,c_j ,a_k ,c_l )$ iff $a_i = a_j$. \\ 
In this case we say that the $n$-pair $\bar{a}\bar{c}$ labels the tuple $\bar{a}$. 

We now define a class $\mathcal{K}$ of finite $(L\cup L_0)$-structures 
as follows. \\ 
(i) In each structure of $\mathcal{K}$ all the relations $R_n$ are defined on $P$; \\ 
(ii) The $P$-part of any structure from $\mathcal{K}$ 
is isomorphic to a finite substructure of $M_0$; \\ 
(iii) For any $D\in \mathcal{K}$, any $n$ and any $n$-pair from $D$ 
labelling a tuple $\bar{a}$ we have $R_n (\bar{a})$.  

It is obvious that $\mathcal{K}$ is closed under substructures and 
there is a function $f:\omega \rightarrow\omega$ so that for every $n$ 
the number of non-isomorphic sructures of $\mathcal{K}$ of size $n$ 
is bounded by $f(n)$.  
The function $f$ is computable with respect to $\Phi$ and the Ryll-Nardzewski 
function. 

\begin{lem} 
The class $\mathcal{K}$ has the amalgamation (and the joint embedding) property. 
\end{lem} 

{\em Proof.} 
Let $D_1$ and $D_2$ be structures in $\mathcal{K}$ with intersection $C$. 
By induction it is enough to deal with the case where 
$|D_1 \setminus C|$ = $|D_2 \setminus C|$ =1.  
Let $D_i \setminus C =\{ d_i \}$ and $d_1 \not= d_2$. 
There are three cases. \\ 
Case 1. $d_1$ and $d_2$ both satisfy $P$. 
Using that $M_0$ has quantifier elimination we amalgamate the $P$-parts 
of $D_1$ and $D_2$ remaining the $Q$-part and $S$ the same as before. 
By (4) 
there are no new $n$-pairs in the amalgam, for any $n$.  \\ 
Case 2. $d_1$ and $d_2$ both satisfy $Q$.
In this case we just take the free amalgamation 
(without any new tuples in relations). 
By (4) 
there are no new $n$-pairs in the amalgam, for any $n$.  \\ 
Case 3. $d_1$ satisfies $P$ and $d_2$ satisfies $Q$.
In this case we again take the free amalgamation and  by (4) we again 
have that there are no new $n$-pairs in the amalgam, for any $n$.
$\Box$ 

\bigskip 

We now see that by Fra\"{i}ss\'{e}'s theorem, the class $\mathcal{K}$ 
has a universal homogeneous (and $\omega$-categorical ) structure $U$.  
In particular ${\mathcal{K}}/\cong$ coincides with $Age(U)$ 
(= collection of all types of finite substructures of $U$). 

Since $M_0$ is the Fra\"{i}ss\'{e} limit of the class of 
all $P$-parts of structures from $\mathcal{K}$, we see that 
the $P$-part of $U$ is isomorphic to $M_0$. 
Let $N$ be the reduct of $U$ to the language $L$. 
Note that $U$ (thus $M_0$) is definable in $N$. 
Indeed each $R_n$ is definable by the rule: 
$U\models R_n (\bar{a})$ if and only if there is 
an $n$-pair in $N$ which labels $\bar{a}$  
(this follows from the fact that $\mathcal{K}$ contains an $n$-pair for such $\bar{a}$).  

If two tuples $\bar{a}$ and $\bar{b}$ in $M_0$ realise the same type 
in $M_0$ they realise the same quantifier free type in $U$. 
So by quantifier elimination there is an automorphism of $U$ 
(and of $N$) which takes $\bar{a}$ to $\bar{b}$. 
This shows that $M_0$ is a dense relativised reduct of $N$. 

To see the last statement of the theorem consider 
a set $\Phi$ axiomatising $Th(M_0)$. 
Thus the $P$-part of $U$ must satisfy $\Phi$ with respect 
to the relations $R_n$ defined in $N$ as above. 
The remaining axioms of $Th(N)$ (and of $Th(U)$) are just 
the axioms of the universal homogeneous structures of 
the corresponding class satisfying (i) - (iii) as above.  
$\Box$ 

\bigskip 

\begin{Remark} \label{ax} 
The structure $U$ produced in the proof is axiomatised as follows. \\
{\bf Axiomatisation of } $Th(U)$. \\ 
(a) all universal axioms forbidding finite substructures which cannot 
occur in $M_0$; \\ 
(b) all universal axioms stating property (iii) from the proof ; \\ 
(c) all $\exists$-axioms for finite substructures of $M_0$; \\ 
(d) all $\forall\exists$-axioms which realise the property of universal 
homogeneous structures that for any $\mathcal{K}$-structures $A <B$ with 
$A<U$ there is an $A$-embedding of $B$ into $U$. 

Note that for every pair of natural numbers $n$ and $l$ the axioms of (a), (b) 
and (c) with at most $n$ quantifiers in the sublanguage of $L\cup L_0$ of arity $\le l$ 
determine all $n$-element structures from $\mathcal{K}$ in this sublanguage. 
On the other hand by the Ryll-Nardzewski function of $Th(M_0 )$ we can find the arity 
$l_n$ so that all $\mathcal{K}$-embeddings between structures of size $\le n$ 
are determined by their relations of arity $\le l_n$. 
Thus the axioms of (d) with at most $n$ quantifiers can be effectively found 
by the corresponding axioms (a - c) and the Ryll-Nardzewski function. 
Moreover there is an effective procedure which for every natural numbers $n$  
produces all $\forall \exists$-sentences of $Th(U)$ with at most $n$ quantifiers, 
when one takes as the input the axioms of (a) and (c) 
of $U$ with at most $n$ quantifiers. 
\end{Remark}

\section{Finite language and non-G-compact theories}

The following definitions and facts are partially 
taken from \cite{clpz}. 
Let ${\bf C}$ be a monster model of the teory $Th({\bf C})$. 
For $\delta\in\{ 1,2,...,\omega\}$ 
let $E_{L}^{\delta}$ be the finest bounded 
$Aut({\bf C})$-invariant equivalence relation on $\delta$-tuples 
(i.e. the cardinality of the set of equivalence classes is bounded). 
The classes of $E_{L}^{\delta}$ are called Lascar strong types. 
The relation $E^{\delta}_{L}$ can be characterized as follows:
$(\bar{a},\bar{b})\in E^{\delta}_{L}$ if there are $\delta$-tuples
$\bar{a}_{0}(=\bar{a}),\bar{a}_1 ,...,\bar{a}_{n}(=\bar{b})$ such that
each pair $\bar{a}_{i},\bar{a}_{i+1}$, $0\le i<n$, extends to
an infinite indiscernible sequence. 
In this case denote by $d(\bar{a},\bar{b})$ the minimal $n$ such 
that some $\bar{a}_0 (=\bar{a} ), \bar{a}_1 ,..., \bar{a}_n (=\bar{b} )$ 
are as above. 

Let $E^{\delta}_{KP}$ be the finest bounded type-definable
equivalence relation on $\delta$-tuples. 
Classes of this equivalence relation are called KP-strong types. 
The theory $Th({\bf C})$ is called {\em G-compact} if 
$E^{\delta}_{L}= E^{\delta}_{KP}$ 
for all $\delta$. 
The first example of a non-G-compact theory was found in \cite{clpz}. 
The first example of an $\omega$-categorical non-G-compact 
theory was found by the author in \cite{ivanov}. 
The following proposition is a straightforward application of Theorem \ref{Udi}. 

\begin{prop} 
There is a countably categorical structure $N$ in 
a finite language such that $Th(N)$ is not G-compact. 
\end{prop} 

{\em Proof.} 
Let $L$ be defined as in the proof of Theorem \ref{Udi}. 
Corollary 1.9(2) of \cite{newelski} states that G-compactness is 
equivalent to existence of finite bound on the diameters of Lascar 
strong types. 
Let $M_0$ be an $\omega$-categorical structure 
which is not G-compact, see \cite{ivanov}. 
In \cite{ivanov} for every $n$ a pair $\bar{a}_n$, 
$\bar{b}_n$ of finite tuples of the same Lascar strong type 
is explicitely found so that $d(\bar{a}_n ,\bar{b}_n )>n$. 

Let $N$ be an $L$-structure, so that $M_0$ is a dense 
relativised reduct in $N$ defined by $P$. 
Then $Th(N)$ is not G-compact. 
Indeed for every $n$, the pair $\bar{a}_n, \bar{b}_n$ is of the same 
Lascar strong type and $d(\bar{a}_n ,\bar{b}_n )>n$ with 
respect to the theory of $N$. 
To see this notice that if in  
$\bar{c}_0 (=\bar{a}_n ), \bar{c}_1 ,..., \bar{c}_m (=\bar{b}_n )$  
each $\bar{c}_i ,\bar{c}_{i+1}$ extends to an indiscernible sequence 
in $Th(M_0 )$, then this still holds in $Th(N)$ by density 
of the image of $Aut(N)$ in $Aut(M_0 )$. 
On the other hand since $Aut(N)\le Aut(M_0 )$ on $P(M)$,  
we cannot find in $N$ such a sequence with $m\le n$.  
$\Box$

\section{Finite language and Ehrenfeucht theories}

In this section we consider the situation where $M_0$ 
is obtained by an $\omega$-sequence of $\omega$-categorical 
expansions. 
We will see that under some natural assumptions the 
construction of Section 1 still works in this 
situation. 
Using this we will prove that there is a finite language such 
that the indexes of Ehrenfeucht theories with exactly 
three countable models form a $\Pi^1_1$-hard set. 

Let $M_0$ be a countable structure of a 1-sorted, 
relational language  $L_0$ = $\{ R_1 ,R_2 ,..., R_n , ...\}$.  
Suppose $L_0 = \bigcup_{i>0} L_i$, where for each $i>0$,  
$L_i = \{ R_1 ,..., R_{l_i}\}$ and the $L_i$-reduct of $M_0$ 
admits quantifier elimination (and thus $\omega$-categorical). 
We may assume that the arity of $R_n$ is not greater than $n$. 
Admitting $R_n$ with repeated coordinats, 
we may also assume that for all $m<n$ the arity of $R_m$ is not 
greater than the arity of $R_n$ and the arity of $R_{l_i}$ 
is less than the arity of $R_{l_{i}+1}$. 

We now admit that $M_0$ is not $\omega$-categorical. 
On the other hand the theory of $M_0$ can be axiomatised 
as follows. 
For each $i$ consider the $L_i$-reduct of $M_0$ and 
its age $Age(M_0 |L_i )$.  
Then this reduct is axiomatised by the standard axioms 
of a universal homogeneous structure (i.e. the versions of 
(a),(c),(d) from Remark \ref{ax} with respect to $Age(M_0 |L_i )$). 
The collection of all systems of axioms of this kind 
gives an axiomatisation of $Th(M_0 )$. 

Applying the proof of Theorem \ref{Udi} we associate 
to each $L_i$-reduct of $M_0$, a class ${\mathcal{K}}_i$ of 
$(L\cup L_i )$-structures obtained by conditions (i)-(iii) 
from this proof. 
Since the $L_i$-reduct of $M_0$ has quantifier elimination, 
repeating the argument of Theorem \ref{Udi} we obtain an  
$\omega$-categorical $(L\cup L_i )$-structure $U_i$ and 
the corresponding $L$-reduct $N_i$ (since the language is finite, 
we do not need the assumption that each $R_i$ describes a type). 
Notice that the construction forbids $n$-pairs 
for $R_n$ of arity greater than the arity of $L_i$. 

\begin{lem} \label{lemexp}
(1) For any $i<j$ the structures $U_{i}$ and $U_{j}$ satisfy 
the same axioms of the form (a) - (d) of Remark \ref{ax} where 
the language of the $P$-part is restricted to $L_{i}$ and the number 
of variables of the $Q$-part is bounded by the arity of $L_i$. 

(2) The corresponding structures $N_{i}$ and $N_{j}$ satisfy 
the same sentences which are obtained by rewriting 
of the axioms of statement (1) as $L$-sentences 
(using the interpretation of $U_i$ in $N_i$). 
\end{lem} 

{\em Proof.} 
Let $m$ be the arity of $L_i$. 
To see statement (1) let us prove that the classes ${\mathcal{K}}_i$ 
and ${\mathcal{K}}_j$ consist of the same $(L\cup L_i )$-structures 
among those with the $Q$-part of size $\le m$. 
The direction $j\rightarrow i$ is clear: the $(L\cup L_i )$-reduct 
of an $(L\cup L_j )$-structure of this form obviously satisfies 
the requirements (i) - (iii) corresponding to ${\mathcal{K}}_j$ 
(and to ${\mathcal{K}}_i$ too). 
To see the direction $i\rightarrow j$ note that the assumption 
that the size of the $Q$-part is not geater than $m$ implies 
that such a structure from ${\mathcal{K}}_i$ has an expansion 
to an $(L\cup L_j )$-structure from ${\mathcal{K}}_j$. 

Now the case of axioms of the form (a),(b),(c) is easy. 
Consider case (d). 
Since the $L_i$-reduct of $M_0$ admits elimination of quantifiers, 
for any finite $L_j$-substructure $A<M_0$ and any embedding of 
the $L_i$-reduct of $A$ into any $B\in Age(M_0 |L_i )$ there 
is an $L_j$-substructure of $M_0$ containing $A$ with the $L_i$-reduct 
isomorphic to $B$.  
This obviously implies that for any substructure $A'<U_j$ without 
$n$-pairs for arities greater than $arity(L_i )$, any embedding of 
the $(L\cup L_i )$-reduct of $A'$ into any $B'\in {\mathcal{K}}_i$ 
can be realised as a substructure of $U_j$ containing $A'$ with 
the $(L\cup L_i )$-reduct isomorphic to $B'$. 
This proves (1). 

Statement (2) follows from statement (1). 
$\Box$ 

\bigskip 

We now additionally assume that $M_0$ is a generic structure 
with respect to the class ${\mathcal{K}}_0$ of all finite 
$L_0$-substructures of $M_0$. 
This means that ${\mathcal{K}}_0$ has the joint embedding 
and amalgamation properties (JEP and AP), 
$({\mathcal{K}}_0 /\cong )=Age(M_0 )$ and $M_0$ is a countable 
union of an increasing chain of structures from ${\mathcal{K}}_0$ 
so that any isomorphism between finite substructures extends 
to an automorphism of $M_0$. 

Let $\mathcal{K}$ be the class of all finite $(L_0 \cup L)$-structures 
satisfying the conditions (i)-(iii) with respect to ${\mathcal{K}}_0$.   
In particular it obviously contains only countably many 
isomorphism types and the class ${\mathcal{K}}_0$ appears 
as the class of all $P$-parts of ${\mathcal{K}}$.  
Applying the proof of Theorem \ref{Udi} we see that ${\mathcal{K}}$ 
is closed under substructures and has the joint embedding and 
amalgamation properties.  
By Theorem 1.5 of \cite{kuelas} the class ${\mathcal{K}}$ 
has a unique (up to isomorphism) generic structure 
(i.e. a structure which is a countable union of an 
increasing chain of structures from ${\mathcal{K}}$ 
and satisfies axioms (a) - (d) of Remark \ref{ax}). 
Note that this structure can be non-$\omega$-categorical.

\begin{lem} \label{Robinson1}
Under the circumstances of this section let $U$ be a generic 
$(L\cup L_0 )$-structure for ${\mathcal{K}}$ as above.  
\parskip0pt  

Then the $P$-part of $U$ is isomorphic to $M_0$. 
The structure $M_0$ is a dense relativised reduct of $U$. 
\end{lem}

{\em Proof.} 
The firts statement is obvious. 
The second statement is an application of back-and-forth. 
$\Box$ 

\bigskip 

It is worth noting here that for every $m$ the amalgamation of 
Theorem \ref{Udi} preserves the subclass of $\mathcal{K}$ 
consisting of structures without $n$-pairs for arities greater than $m$ 
(for example structures with the size of the $Q$-part less than $m+1$).  
If $m$ is the arity of the language $L_i$ then 
$(L_i \cup L)$-reducts of these structures form the Fra\"{i}ss\'{e} 
class corresponding to the universal homogeneous structure $U_i$. 

\begin{prop} \label{Robinson} 
(1) All axioms of $U_i$ of the form (a),(c),(d) of 
Remark \ref{ax} also hold in $U$. \\ 
(2) The theory $Th(U)$ is model complete and is axiomatised 
by axioms of the form (b) of Remark \ref{ax} together with 
the union of all axioms of the form (a),(c),(d) for all $Th(U_i)$. \\  
(3) For any axiom $\phi$ of $Th(U)$ of the form (a)-(d) as in (2)  
there is a number $i$ so that $\phi$ holds in all $U_j$ for $j>i$. 
\end{prop} 

{\em Proof.} 
(1) The case of axioms of the form (a),(c) is easy. 
Consider case (d). 
Since the $L_i$-reduct of $M_0$ admits elimination of quantifiers, 
for any substructure $A<M_0$ and any embedding of the $L_i$-reduct 
of $A$ into any $B\in Age(M_0 |L_i )$ there is a substructure of 
$M_0$ containing $A$ with the $L_i$-reduct isomorphic to $B$.  
This obviously implies that for any substructure $A'<U$ without 
$n$-pairs of arity greater than $arity(L_i )$, any embedding of 
the $(L\cup L_i )$-reduct of $A'$ into any $B'\in {\mathcal{K}}_i$ 
can be realised as a substructure of $U$ containing $A'$ with 
the $(L\cup L_i )$-reduct isomorphic to $B'$. 
This proves (1). 

(2) Let $U'$ and $U''$ satisfy axioms as in the formulation of (2). 
Then obviously the $(L_i \cup L)$-reducts of 
$U'$ and $U''$ satisfy the axioms of $Th(U_i )$ as in statement (1). 
In particular $P(U')\cong P(U'' )$ in each $L_i$. 
Moreover if $U' <U''$, then by axioms (d) one can easily verify 
that this embedding is $\forall$-elementary.  
Thus $U'$ is an elementary substructure of $U''$ by 
a theorem od Robinson. 
It is also clear that $U$ is embeddable into any structure 
satisfying axioms as in (2). 

(3) By Lemma \ref{lemexp} we see that for every sentence 
$\theta \in Th(U)$ of the form (a) - (d) of (2) there is 
a number $i$ such that for all $j>i$, $\theta$ holds in $U_j$. 
$\Box$ 

\bigskip 

Some typical examples of Ehrenfeucht theories 
(i.e. with finitely many countable models) are build 
by the method of this section:  
the theory of all expansions of $(\mathbb{Q},<)$ by infinite 
discrete sequences $c_1 < c_2 < ...<c_n <...$, is Ehrenfeucht 
and can be easily presented in an appropriate $L_0$ as above.

\begin{prop} \label{Ehrenfeucht1}
Under the circumstances of this section assume that $M_0$ 
is a generic structure with respect to the class $\mathcal{K}_0$ 
of all finite substructures of $M_0$.  
Assume that $Th(M_0 )$ is an Ehrenfeucht theory. 
Let $U$ be a generic $(L\cup L_0 )$-structure for ${\mathcal{K}}$ as above.  

Then $Th(U)$ is also Ehrenfeucht. 
\end{prop} 

{\em Proof.} 
Let $U',U''$ be countable models of $Th(U)$. 
Assume that the $P$-parts of $U'$ and $U''$ 
(say $M'$ and $M''$) are isomorphic. 
Identifying them let us show that $U'$ is isomorphic to $U''$. 
For this we fix a sequence of finite substructures 
$A_1 <A_2 <...<A_i <...$ so that $M' =\bigcup A_i$.  
Having enumerations of the $Q$-parts of $U'$ and $U''$ 
we build by back-and-forth, sequences $B'_1 <B'_2 <...<B'_i <...$ and 
$B''_1 <B''_2 <...<B''_i <...$ with $B'_i >A_i <B''_i$, $U' =\bigcup B'_i$   
and $U'' =\bigcup B''_i$ so that $B'_i$ is isomorphic to $B''_i$ 
over $A_i$. 
By Proposition \ref{Robinson}(2) using the fact that 
$U',U''\models Th(U)$ we see that such sequences exist.  
$\Box$ 

\bigskip 

We now prove that there is a finite language $L$ such that the set of 
Ehrenfeucht $L$-theories with exactly three models is $\Pi^1_1$-hard.

\begin{thm} \label{Ehrenfeucht2} 
There is a finite language $L$ such that for every  
$B \in \Pi^1_1$ there is a Turing reduction of $B$ to 
the set $3Mod_L$  of all indexes of decidable Ehrenfeucht 
$L$-theories with exactly three countable models. 
\end{thm}  

{\em Proof.} 
Let $L$ be the language defined in Section 1. 
We use the idea of Section 4 of \cite{LS}. 
In particular we can reduce the theorem to the case when 
$B$ coincides with the index set $NoPath$ of the property  
of being a computable tree $\subseteq \omega^{\omega}$ having 
no infinite path. 
The Turing reduction of this set to $3Mod_L$ which will be built below, 
is a composition of the procedure described in 
\cite{reed} and \cite{LS}, and the construction of this section.  
The former one is as follows. 
Having an index $e$ of a computable tree $Tr_e \subset \omega^{\omega}$,  
R.Reed defines a complete decidable theory $T_e$ of the language 
$$ 
\langle \wedge, <_{L}, \le_{H}, E^{\eta}_{\xi}, L^{\eta}_{\xi}, H_{\eta}, A_{\eta}, B_{\eta}, c_{\eta}    
\mbox{ ( } \eta , \xi\in Tr_e \mbox{ ) }\rangle , 
$$ 
where $\wedge$ is the function of the greatest lower bound of a tree, 
$<_{L}$ is a Kleene-Brouwer ordering of this tree and $\le_H$ is 
a binary relation  measuring 'heights' of nodes.  
Constants $c_{\eta}$, $\eta \in Tr_e$, define embeddings of $Tr_e$ 
into models of $T_e$.       
The remaining relations are binary. 

For each natural $n$ define $T_e |_n$ to be the restriction of $T_e$ 
to the sublanguage corresponding to the indexes from the finite subtree 
$Tr_e \cap n^{<n}$. 
The proof of Lemma 9 from \cite{reed} shows that $T_e |_n$ admits 
effective quantifier elimination.  
Lemma 6 of \cite{reed} asserts that every quantifier-free formula 
of $T_e |_n$ is equivalent to a Boolean combination of atomic formulas 
of the following form:  
$$ 
u\wedge v = w\wedge z \mbox{ , } u<_{L} w \mbox{ , } u\wedge v \le_{H} w\wedge z \mbox{ , } 
E^{\eta}_{\xi}(u,w) \mbox { , } 
$$ 
$$
L^{\eta}_{\xi}(u,w) \mbox{ , } H_{\eta}(u\wedge v,w) \mbox{ , } A_{\eta}(u\wedge v,w) \mbox{ , } 
$$ 
where $u,v,w,z$ is either a variable or a constant in  $T_e |_n$. 
By Lemma 8 of \cite{reed} the corresponding Boolean combination 
can be found effectively.  
This implies that replacing the function $\wedge$ by the first, 
third, sixth and seventh relations of the list above we transform 
the language of each $T_e$ into an equivalent relational language. 
In particular we have that each $T_e |_n$ is $\omega$-categorical.  

Note that extending the set of relations we can eliminate 
constants $c_{\eta}$ from our language. 
Admitting empty relations we may assume that all $T_e$ have 
the same language (where $\omega^{<\omega}$ is the set of indexes).   
Admitting repeated coordinates we may assume that this language 
$L_0 = \{ R_1, ..., R_i ,...\}$ satisfies the assumptions of 
the beginning of the section and each sublanguage $L_n$ of 
the presentation $L_0 = \bigcup_{i>0} L_i$ corresponds to $T_e |_n$. 

We now apply Lemma \ref{lemexp} to all $T_e |_n$. 
Since each $T_e |_n$ is computably axiomatisable uniformly 
in $e$ and $n$, we obtain an effective enumeration 
of computable axiomatisations of $L$-expansions of all 
$T_e |_n$ (with $T_e |_n$ on the $P$-part).   
For each $e$ taking the axioms which hold in almost all 
$L$-expansions of $T_e |_n$ we obtain by Lemma \ref{lemexp}(1)  
a computable axiomatisation of a theory of $L$-expansions of $T_e$. 

When $T_e$ is an Ehrenfeucht theory with exactly 
three models (i.e. $e\in NoPath$), the prime model 
of $T_e$ is generic with respect to its age.  
Applying Proposition \ref{Robinson} to 
$T_e$ and all $T_e |_n$ we obtain a generic 
$(L_0 \cup L)$-structure such that its theory is 
computably axiomatised as above.   
This theory has exactly three countable models 
by Proposition \ref{Ehrenfeucht1}. 

When we take $L$-reducts of all $T_e$ and 
the corresponding computable axiomatisations 
we obtain a computable enumeration of $L$-theories 
which gives the reduction as in the formulation of the theorem. 
$\Box$ 

\bigskip

\begin{Remark} 
In the proof above we used Proposition \ref{Robinson} 
in order to obtain a complete $L$-expansions of Ehrenfeuch $T_e$'s.  
We cannot apply it in the case when $T_e$ does not have 
an appropriate generic model, for example when the 
corresponding $Tr_e$ has continuum many paths. 
Nevertheless the author hopes that 
the proof can be modified so that the 
reduction as above also shows that the set of 
all $L$-theories with continuum many models 
is $\Sigma^1_1$-hard. 
In the case of infinite languages this is shown 
in Section 4 of \cite{LS}. 
\end{Remark}

\section{Coding $\omega$-categorical theories} 

The main theorem of this section improves the corresponding result of 
\cite{LS} (where the authors do not demand that the language is finite). 
It is worth noting that the author together with Barbara Majcher-Iwanow 
have found some other improvements in \cite{IMI}.   

\begin{thm} 
There is a finite language $L$ such that the property of 
$\omega$-categoricity distinguishes a $\Pi^{0}_3$-complete 
subset of the set of all decidable complete $L$-theories. 
\end{thm}  

{\em Proof.} 
In the formulation of the theorem $L$ is the 
language defined in Section 1. 
It is shown in \cite{LS} that the property of 
$\omega$-categoricity is $\Pi^{0}_3$. 
The proof of $\Pi^{0}_3$-completeness in the case of $L$ is based 
on Theorem \ref{Udi}, Section 3 and  the idea of Section 2 of \cite{LS}. 
The latter one will be presented in some special form, 
the result of a fusion with some ideas from \cite{P}.  

Let us fix the standard enumeration $p_n$ of prime numbers and 
a G\"{o}del 1-1-enumeration of the set of pairs $\langle i,j\rangle$. 
Let $a(x)$ be a computable increasing function 
from $\omega$ to $\omega \setminus \{ 0,1,2\}$ so that 
if natural numbers $x_1 <x_2$ enumerate pairs 
$\langle i_1 ,j_1 \rangle$ and $\langle i_2 ,j_2 \rangle$ 
then $p_{i_1}a(x_1) <p_{i_2}a(x_2 )$. 

Let $L_E$ consist of $2p_n$-ary relational symbols $E_n$, 
$n\in\omega$, and $T_E$ be the 
$\forall\exists$-theory of the universal homogeneous 
structure of the universal theory saying that each $E_n$ is 
an equivalence relation on the set of $p_n$-tuples 
which does not depend on the order of tuples and 
such that all $p_n$-tuples with at least one repeated 
coordinate lie in one isolated $E_n$-class (Remark 4.2.1 in \cite{P}).  
It is worth mentioning here that the joint embedding 
property and the amalgamation property are easily 
verified by an appropriate version of free amalgamation 
(modulo transitivity of $E_n$-s).  
Note also that $T_E$ is $\omega$-categorical and decidable. 
\parskip0pt 

We now define an auxiliary language $L_{ESP}$. 
We firstly extend $L_E$ by countably 
many sorts $S_n$, $n\in\omega$. 
Start with a countable model $M_E \models T_E$ 
and take the expansion of $M_E$ to the language 
$L_E \cup \{ S_1, ...,S_n ,...\} \cup \{ \pi_1 ,..., \pi_n ,...\}$, 
where each $S_n$ is interpreted by the non-diagonal elements 
of $M^{p_n} /E_n$ and $\pi_n$ by the corresponding projection. 
To define $L_{ESP}$ we extend 
$L_E \cup \{ S_1 ,...,S_n ,...\} \cup \{ \pi_1 ,...,\pi_n ,...\}$ 
by an $\omega$-sequence of relations $P_{m}$, $m\in\omega$, with 
the following properties. 
If $m$ is the G\"{o}del number of the pair $\langle n,i\rangle$ 
then we interpret $P_m$ as a subset of the diagonal of $S^{a(m)}_n$. 
Let $T_{ESP}$ be the $L_{ESP}$-theory axiomatized by $T_E$ 
together with the natural axioms for all $\pi_n$ and $P_m$ as above.  

Having a structure $M\models T_{ESP}$ (which is an expansion of $M_E$) 
we now build another expansion $M^*$ of $M_E$ (in the 1-sorted language). 
For each relational symbol $P_{m}$ of the sort 
$S^{a(m)}_{n}$ we add a new relational symbol $P^*_{m}$ on 
$M^{a(m) p_n}_E$ interpreted in the following way: 
$$ 
M^* \models P^*_{m} (\bar{a}_1 ,...,\bar{a}_{a(m)} ) \Leftrightarrow 
M\models P_{m} (\pi_{n}(\bar{a}_1 ),...,\pi_{n}(\bar{a}_{a(m)} )). 
$$  
It is clear that $M^*$ and $M$ are bi-interpretable. 

By $T^*_{ESP}$ we denote the theory of all $M^*$ with $M\models T_{ESP}$. 
Let $L_0$ be the corresponding language. 
Then $M_E$ is the $L_E$-reduct of any countable $M^*\models T^*_{ESP}$. 
It is clear that $T^*_{ESP}$ is axiomatized by the 
$\forall\exists$-axioms of $T_E$, $\forall$-axioms 
of $E_n$-invariantness of all $P^*_{m}$ and $\forall$-axioms 
that every $P_m$ is a subset of an appropriate diagonal.  
Moreover for every natural $l$ we have $\le 1$ relations 
of arity $l$ in $L_0$ and the function of arities of $P^*_m$ 
is increasing. 
Admitting empty relations (say $R_j$) we may think that for every natural 
number $l>0$ the language $L_0$ contains exactly one relation of arity $l$. 
In particular $L_0$ satisfies basic requirments 
on $L_0$ from Section 3.  
We present $L_0$ as the union of a sequence of finite languages 
$L_1 \subset L_2 \subset ... \subset L_m \subset ...$ of 
arities $l_1 <l_2 <...<l_m <...$ where $L_m$ consists of all 
relations of arity $\le p_n a(m)$ ($=l_m$) with $n$ 
to be the first coordinate of the pair enumerated by $m$. 
Note that when $m$ codes a pair $\langle n ,j\rangle$ 
the relation $E_n$ is also in $L_m$. 
\parskip0pt 

For every $m\in \{ 1 ,...,i ,..., \omega\}$ and a finite 
set $D$ of indexes of relations $P^*_i$ of arity $\le l_m$ 
we consider the class $\mathcal{K}_{D}$ of all finite substructures 
of models of $T^*_{ESP}$ satisfying the property that all 
$P^*_{i}$ with $i\not\in D$, are empty. 
It is clear that for any natural number $k$ the number of 
structures of $\mathcal{K}_{D}$ of size $k$ is finite. 
We will also denote  $\mathcal{K}_{\omega ,D}:=\mathcal{K}_{D}$. 
When $m <\omega$ we define $\mathcal{K}_{m ,D}$ as the class 
of all reducts of $\mathcal{K}_{D}$ to the sublanguage $L_m$.  

By an appropriate version of free amalgamation we see 
that $\mathcal{K}_{m,D}$ has the joint embedding property 
and the amalgamation property. 
Let $M_{m,D}$ be the corresponding universal homogeneous 
structure and let $T^*_{m,D}$ be the theory of $M_{m,D}$. 
It follows from $T^*_{\omega ,D}$ that $T^*_{ESP} \subset T^*_{\omega ,D}$ 
and for every $n$ the family of all $P_{i}$ 
\footnote{$L_{ESP}$-predicates corresponding to $P^*_i$}
, with $i\in D$ coding some $\langle n ,j\rangle$, 
freely generates a Boolean algebra of infinite subsets 
of the sort $S_n$  
(we may interpret such $P_i$ as a unary predicate on $S_n$). 
\parskip0pt 
    
By the definition of the class $\mathcal{K}_{D}$ we see that 
for any $t< m$ and any two finite sets $D'$ and $D''$ satisfying  
\begin{quote} 
$D' \cap \{ 0,...,l_t \} = D'' \cap \{ 0,...,l_t \}$ 
\end{quote} 
the reducts of $M_{m,D'}$ and $M_{m,D''}$ to 
$L_t$ are isomorphic. 

Let us apply the construction of Theorem \ref{Udi} to $M_{m,D}$. 
Then we obtain the $(L_m \cup L )$-structure $U_{m,D}$ and 
the corresponding $L$-reduct $N_{m,D}$, where $L$ is 
the language as in Theorem \ref{Udi}. 
It follows from the proof of that theorem that in 
the situation of the previous paragraph the structures 
$U_{m,D'}$ and $U_{m,D''}$ satisfy the same axioms of 
the form (a) - (d) of Remark \ref{ax}, where the language 
of the $P$-part is restricted to $L_t$ and the number 
of variables of the $Q$-part is bounded by $l_t$. 
When we rewrite these axioms as $L$-sentences 
(using the corresponding definition of the relations 
of $L_m$) we obtain that $N_{m,D'}$ and $N_{m,D''}$ 
satisfy the same axioms of this kind. 

Let $\varphi (x,y)$ be a universal computable function, 
i.e. $\varphi (e,x)=\varphi_e (x)$. 
Find a computable function $\rho$ (with $Dom(\rho )=\omega$) 
enumerating $Dom (\varphi (\varphi (y,z),x))$,   
i.e. the set of all triples $\langle e,n,x\rangle$ 
with $x\in W_{\varphi_{e}(n)}$. 

For any natural $e,s$ we define a finite set $D^s_e$ of 
codes $m\le l_s$ of all pairs $\langle n ,k\rangle$ such that 
$$ 
(\exists x)(\rho (k)=\langle e,n ,x\rangle \wedge 
(\forall k'<k)(\rho (k')\not=\langle e, n ,x\rangle )).
$$ 
Let $T_e$ and $T^*_e$ be the $L_{ESP}$-theory and 
the corresponding 1-sorted version (containing $T^*_{ESP}$) 
such that for all natural $s$   
the reduct of $T^*_e$ to $L_s$ 
coincides with the corresponding reduct of $T^*_{s,D^s_e }$. 
Since for any $s<t$ we have $D^{t}_e \cap \{ 0,...,l_s\} = D^s_e$, 
the definition of $T_e$ and $T^*_e$ is correct. 
It is clear that both $T_e$ and $T^*_e$ are axiomatisable 
by computable sets of axioms uniformly in $e$. 
Since for each $s$ the reduct of $T^*_e$  
as above is $\omega$-categorical, the theories $T_e$ 
and $T^*_e$ are complete. 
Thus $T_e$ and the corresponding theory $T^*_e$ are   
decidable uniformly in $e$.  
It is worth noting that for each $m$ the $L_m$-reduct 
of $T^*_e$ admits elimination of quantifiers 
(it is of the form $T^*_{m,D}$ as above). 
Moreover, the class 
$\bigcup_{l} \mathcal{K}_{\omega,D^l_e}$ 
considered as a class of $L_0$-structures 
where almost all $P^*_m$ are empty, 
is a countable class with JEP and AP. 
It is clear that $T^*_e$ is 
the theory of the corresponding universal
homogeneous structure $M^*_e$. 

Applying Proposition \ref{Robinson} to 
$M^*_e$ and all $M_{l,D^l_e}$ we obtain the 
$(L_0 \cup L)$-structures $U_e$ and their 
approximations $U_{l,D^l_e}$ (and $N_{l,D^l_e}$),  
which for $l\rightarrow \infty$ give a computable 
axiomatisation of the complete $L$-theory $T^L_e$
of the corresponding $L$-reducts $N_e$. 
By Remark \ref{ax} applied to all $U_{l,D^l_e}$ 
(with decidable theories),  this axiomatisation 
(the corresponding decidability of $T^L_e$) can be 
found by an effective uniform in $e$ procedure.  

We now fix a G\"{o}del coding of the language $L$, and identify 
decidable complete $L$-theories with computable functions from 
$\{ sgn(\varphi_e (x)): e\in\omega\}$ realising the corresponding 
characteristic functions (by $sgn(x)$ we denote the function 
which is equal to $1$ for all non-zero numbers and $sgn(0)=0$).  
We want to prove that the set of 
all natural numbers $e$ satisfying the relation 
 
"$sgn(\varphi_e (x))$ codes a decidable $\omega$-categorical  theory"  \\  
is $\Pi^{0}_3$-complete. 

Fix a Turing machine $\kappa (x,y)$ which decides when for a pair 
$d,e$ the number $d$ codes a sentence which belongs to $T^L_e$ 
(in this case $\kappa (d,e)=1$). 
The following procedure defines a computable function $\xi (z)$ and 
a computably enumerable set $Z$. 
At step $e$ we take the Turing machine for $sgn(\varphi_e (x))$ 
and check if any replacement of some parameter $e'$  
in that program by a variable $y$ makes it the Turing machine 
$\kappa (x,y)$. 
If this happens we put $e$ into $Z$ and define $e':=\xi (e)$. 
As a result we obtain a computably enumerable set $Z$ and 
a computable function $\xi$ with $Dom(\xi )\supset Z$ and 
$Rng (\xi )=\omega$ such that for every $e\in Z$ the function 
$sgn (\varphi_{e} (x))$ is computed by the machine 
$\kappa (x,\xi (e))$ (for $T^L_{\xi (e)}$).   

By Ryll-Nardzewski's theorem the $L_{ESP}$-theory $T_{\xi (e)}$ is 
$\omega$-categorical if and only if all $W_{\varphi_{\xi (e)}(n)}$ 
are finite (i.e. the set of 1-types 
(pairwise non-equivalent Boolean combinations of $P_{m}$) 
of each $S_n$ is finite). 
If we consider the corresponding $T^L_{f(e)}$, then this property remains true. 
\parskip0pt

Since for any Turing machine computing $\varphi_{e'} (x)$ 
we can effectively find a Turing machine deciding $T^L_{e'}$ 
(i.e. in fact we can find $sgn(\varphi_e (x))$ with $\xi (e)=e'$), 
we see that the $\Pi^{0}_{3}$-set 
$\{ e': \forall n (W_{\varphi_{e'} (n)}$ is finite)$\}$ 
is reducible to 
$\{ e: sgn(\varphi_e (x))$ codes an $\omega$-categorical $L$-theory$\}$.  
Since the former one is $\Pi^{0}_{3}$-complete (see \cite{LS}) 
we have the theorem. 
$\Box$

\bigskip 

INSTITUTE OF MATHEMATICS, UNIVERSITY OF WROC{\L}AW, \parskip0pt 

pl.GRUNWALDZKI 2/4, 50-384 WROC{\L}AW, POLAND \parskip0pt

E-mail: ivanov@math.uni.wroc.pl

\end{document}